\documentclass[12pt, reqno]{amsart}
\usepackage{amsmath, amsthm, amscd, amsfonts, amssymb, graphicx, color}
\usepackage[bookmarksnumbered, colorlinks, plainpages]{hyperref}

\newtheorem{theorem}{Theorem}[section]
\newtheorem{lemma}[theorem]{Lemma}
\newtheorem{proposition}[theorem]{Proposition}
\newtheorem{corollary}[theorem]{Corollary}
\theoremstyle{definition}
\newtheorem{definition}[theorem]{Definition}
\newtheorem{example}[theorem]{Example}

\theoremstyle{remark}

\numberwithin{equation}{section}

\usepackage{pdfsync}

\begin{document}
\setcounter{page}{1}

\title[A generalized Birkhoff-James orthogonality and norm parallelism]{A generalized Birkhoff-James orthogonality and norm parallelism in unital $C^*$-algebras and their characterizations}

\author[H. S. Jalali Ghamsari]{Hooriye Sadat Jalali Ghamsari$^{1}$}
\address{$^1$ H. S. Jalali Ghamsari : Department of Pure Mathematics, Faculty of Mathematical
Sciences, University of Kashan, Kashan, 87317-53153, Iran.}
\email{jalali.hooriyesadat@gmail.com}
\author[M. Dehghani]{Mahdi Dehghani$^{2}$}
\address{$^2$ M. Dehghani (corresponding author): Department of Mathematical Sciences, Yazd University, Yazd, Iran.}
\email{e.g.mahdi@gmail.com, m.dehghani@yazd.ac.ir}


%%%%%%%%%%%%%%%%%%%%%%%%%%%%%%%%%%%%%%%%%%%%%%%%%%%%%
\begin{abstract}

Let $\mathcal{A}$ be a unital $C^*$-algebra and let $a\in\mathcal{A}$ be a positive and invertible element. Suppose that $\mathcal{S}(\mathcal{A})$ is the set of all states on $\mathcal{\mathcal{A}}$ and
let
\[ \mathcal{S}_a (\mathcal{A})=\{ \dfrac{f}{f(a)} \, : \, f \in \mathcal{S}(\mathcal{A}), \, f(a)\neq 0\}.\]
We introduce a family of generalized norms, called $(a,\lambda)$-norms, on $\mathcal{A}$ defined by
\[\|x\|_{a,\lambda}:=\sup\{\sqrt{\lambda\varphi(x^*ax)+(1-\lambda)|\varphi(ax)|^2} : \varphi\in\mathcal{S}_a(\mathcal{A})\}\quad(\lambda\in [0,1]).\]
These family of norms generalizes
the recently introduced $a$-operator norm, $\|\cdot\|_a$ and $a$-numerical radius norm, $v_a(\cdot)$ in unital $C^*$-algebras.

The notions of Birkhoff-James orthogonality and norm-parallelism with respect to $\|\cdot\|_{a,\lambda}$, which is called, $(a,\lambda)$-Birkhoff-James orthogonality and $(a,\lambda)$-norm parallelism in $\mathcal{A}$, respectively, are introduced and investigated.
Characterizations of $(a,\lambda)$-norm parallelism and $(a,\lambda)$-Birkhoff-James orthogonality in terms of the elements of $\mathcal{S}_a(\mathcal{A})$ are obtained. In particular, the relationship between these new concepts are described. Our results extend and cover some known results in this area.
 
\end{abstract} \maketitle

%%%%%%%%%%%%%%%%%%%%%%%%%%%%%%%%%%%%%%%%%%%%%%%%%
\section{Introduction and preliminaries}

Let $ \mathcal{A} $ be a unital $ C^{*} $-algebra with  unit  $1_{\mathcal{A}}$. We denote by $\mathcal{A}^{\prime}$  the topological dual space of $\mathcal{A}$. The adjoint of any element $ x \in \mathcal{A} $ is denoted by $ x^* $. Also, for
an element $x\in\mathcal{A}$ we denote by ${\rm Re}(x)=\dfrac{1}{2}(x+x^*)$ and ${\rm Im}(x)=\frac{1}{2i}(x-x^*)$ the
real and the imaginary part of $x$, respectively. An element $a$ of $\mathcal{A}$ is called positive (written by $a\geq 0$), if $a$ is selfadjoint whose spectrum $\sigma(a)$ is contained in $[0,\infty)$. It is known that if $a\in\mathcal{A}$ is positive, then there exists a unique positive element $b\in\mathcal{A}$ such that $a=b^2$. Such an element $b$ is called the positive square root of $a$ and is denoted by $a^{\frac{1}{2}}$. The symbol $\mathcal{A}^+$  stands for the cone of positive elements in $ \mathcal{A} $. A linear functional $f$ on $\mathcal{A}$ is  called positive if $f(a)\geq 0$ for every positive element $a\in\mathcal{A}$. Given a positive functional $f$ on $\mathcal{A}$, the following well-known version of the
Cauchy-Schwartz inequality holds for every $x,y\in\mathcal{A}$:
\[|f(x^*y)|^2\leq f(x^*x)f(y^*y).\]
It is well-known that a linear functional on $\mathcal{A}$ is positive if and only if $f(1_{\mathcal{A}})=\|f\|$; see \cite[Corollary 3.3.4]{Murphy}. Let $ \mathcal{S}(\mathcal{A}) $ be the set of all states on $ \mathcal{A} $. In fact,
$$ \mathcal{S}(\mathcal{A})  = \{ f\in \mathcal{A}^{'} \, : \, f(1_{\mathcal{A}})= \| f \| =1 \}. $$
Let $ a $ be a nonzero positive element of $ \mathcal{A} $. A generalization of state space of $\mathcal{A}$ was introduced in \cite{AB} as follows:
\begin{align*}
\mathcal{S}_{a}(\mathcal{A}):&= \{ \varphi \in \mathcal{A}^{'} \, : \, \varphi \geq 0, \, \varphi ( a ) = 1\}=\{ \dfrac{f}{f(a)} \, : \, f\in \mathcal{S}(\mathcal{A}), \, f(a) \neq 0 \}.
\end{align*}
Observe that if $a=1_{\mathcal{A}}$, then $\mathcal{S}_a(\mathcal{A})=\mathcal{S}(\mathcal{A})$.
It has been proved in \cite{AB} that $\mathcal{S}_a(\mathcal{A})$ is a nonempty convex and $w^*$-closed subset of $\mathcal{A}^{\prime}$. But, unlike $\mathcal{S}(\mathcal{A})$ , the set $\mathcal{S}_a(\mathcal{A})$
 may not be $w^*$-compact. In fact, according to \cite[Proposition 2.3]{AB}, $\mathcal{S}_a(\mathcal{A})$ is $w^*$-compact if and only if $a$ is invertible.
 For any element $x\in\mathcal{A}$, let
\begin{align*}
\|\cdot\|_a :\mathcal{A}\rightarrow[0,\infty),\quad
 \| x \|_{a} := \sup \{ \sqrt{\varphi (x^* a x )} \, : \varphi \in \mathcal{S}_a (\mathcal{A}) \}.
\end{align*}
Due to \cite{AB}, if $a$ is not invertible, then $\mathcal{S}_a(\mathcal{A})$ is not $w^*$-compact, and so it may happen that
$\|x\|_a=\infty$ for some $x\in\mathcal{A}$; see \cite[Example 3.2]{AB}. Denote by $\mathcal{A}^a:=\{x\in\mathcal{A} : \|x\|_a<\infty\}$. It was shown in \cite{AB} that $\|\cdot\|_a$ is a submultiplicative semi-norm on $\mathcal{A}^a$; i.e., $\|xy\|_a\leq\|x\|_a\|y\|_a$ for all $x,y\in\mathcal{A}^a$. Also, $\|x\|_a=0$ if and only if $ax=0$. In addition, if $ a $ is invertible, then $ \| \cdot \|_a $ is a norm on $ \mathcal{A} $. Consequently,  $\|\cdot\|_{1_{\mathcal{A}}}$ is equal to the $C^*$-norm $\|\cdot\|$ of $\mathcal{A}$.

An element $x^{\sharp}\in\mathcal{A}$ is called  an $a$-adjoint of $x\in\mathcal{A}$ if $ax^{\sharp}=x^*a$. The set of all $a$-adjointable elements of $\mathcal{A}$ is denoted by $\mathcal{A}_a$. Note that $\mathcal{A}_a=\mathcal{A}$ if $\mathcal{A}$ is commutative. Moreover, any element $ x \in \mathcal{A}_a $ can be written as $ x=x_1 + i x_2 $,  where $ x_1 $ and $ x_2 $ are $a$-selfadjoint. But, in general, this decomposition is not unique. In fact, if $ x^{\sharp} $ is an $a$-adjoint of $ x $, then
\[x = \dfrac{x+x^{\sharp}}{2} + i \dfrac{x-x^{\sharp}}{2i}.\]
In \cite[Corollary 4.9]{AB} was proved that if $ x \in \mathcal{A}_a $ and $ x^{\sharp} $ is an $ a $-adjoint of it, then
\begin{align}\label{s0009}
\|x\|_{a}^2 = \|xx^{\sharp}\|_{a} = \|x^{\sharp} x\|_{a} =\|x^{\sharp}\|_{a}^2 .
\end{align}
An element $x\in\mathcal{A}$ is said to be $a$-selfadjoint if $ax$ is hermitian; i.e., $ax = x^*a$. Also,
$x$ is called $a$-positive provided that $ax\in\mathcal{A}^+$.

The algebraic $a$-numerical range of any element $x\in\mathcal{A}$ is defined by
\[V_a(x)=\{\varphi(ax) : \varphi\in\mathcal{S}_a(\mathcal{A})\}.\]
Observe that $V_{1_{\mathcal{A}}}(x)=V(x)=\{f(x) : f\in\mathcal{S}(\mathcal{A})\}$ which is known as algebraic numerical range of $x$. It has been proved in \cite[Theorem 4.7]{AB} that $V_a(x)$ is a nonempty convex subset of complex numbers for all $x\in\mathcal{A}^a$. Unlike the classical algebraic
numerical range $V(x)$, the $a$-numerical range $V_a(x)$ of an element $x\in\mathcal{A}$ may not
be bounded. Note that if $a$ is invertible, then $V_a(x)$ is bounded. The maximum
modulus of $V_a(x)$ is called the $a$-numerical radius $v_a(x)$ of $x$; $ v_a (x)= \{ |z| \, : \, z \in V_a (x) \}$. Also, note that $v_{1_{\mathcal{A}}}(x)=v(x)=\sup\{|z| : z\in V(x)\}$. If $a$ is invertible, then it
has been shown that $v_a(\cdot)$ define a norm on $\mathcal{A}$, which is equivalent to  $\|\cdot\|_a$. In fact, the following inequalities are hold:
\begin{align}\label{Eq0}
\frac{1}{2}\|x\|_a\leq v_a(x)\leq\|x\|_a\quad(\forall x\in\mathcal{A}).
\end{align}
For more information about algebraic $a$-numerical ranges and $a$-numerical radius and their fundamental properties
the reader is referred to \cite{Mabrouk,MZ}.

Let $\mathcal{B}(\mathcal{H})$ be the $C^*$-algebra of all bounded linear operators on Hilbert space $(\mathcal{H},\langle\cdot,\cdot\rangle)$ and let $A\in\mathcal{B}(\mathcal{H})$ be a positive operator. From another point of view, the concepts of $V_a(\cdot)$, $v_a(\cdot)$ and $\|\cdot\|_a$ were introduced as generalizations of  $A$-numerical range, $A$-numerical
radius and $A$-operator semi-norm for operators defined by $W_A(T)=\{\langle ATx,x\rangle : \langle Ax,x\rangle=1\}$, $\|T\|_A=\sup\{\sqrt{\langle ATx,x\rangle} : \langle Ax,x\rangle=1, x\in\mathcal{H}\}$ and  $w_A(T)=\sup\{|z| : z\in W_A(T)\}$, respectively,  for every operator $T\in\mathcal{B}(\mathcal{H})$; see e.g., \cite{AML,AM,Bhunia,Bottazzi,AZ} and their references.

It is well-known that the concept of Birkhoff-James orthogonality provide a good framework for studding the geometry of normed linear spaces and operator spaces; see e.g., \cite{AGKZ,James,Serb} and the references therein.
The concepts of  Birkhoff-James orthogonality and  norm-parallelism with respect to the $C^*$-norm and the algebraic numerical radius norm in $C^*$-algebras and also in Hilbert $C^*$-modules have been extensively discussed and studied in \cite{LA,TB,KZ,Moslehian,Will,PW,Wojcik,WZ,ZM1,ZM2,Z2,ZW}. In particular, characterizations of $A$-Birkhoff-James ($A$-numerical Birkhoff-James) orthogonality and
$A$-norm parallelism ($A$-numerical radius norm parallelism) of bounded linear operators with respect to the semi-norm induced
by a positive operator $A $ acting on a complex Hilbert space have been studied in \cite{Bhunia,Bottazzi,AZ}. Recently, the notion of Birkhof-James orthogonality with respect to  $\|\cdot\|_a$ in unital $ C^* $-algebra $\mathcal{A}$, so called $a$-Birkhoff-James orthogonality, has been investigated in \cite{JD}.  Furthermore, a new semi-norm for $A$-bounded operators on Hilbert spaces which generalizes simultaneously $A$-operator semi-norm and $A$-numerical radius was introduced and studied in \cite{EAZ}.

Given a positive and invertible element $a$ in a unital $C^*$-algebra $\mathcal{A}$, and a number $\lambda\in [0,1]$. Inspiring from the ideas that have been presented in \cite{EAZ}, in this paper we consider a new family of norms $\|\cdot\|_{a,\lambda}:\mathcal{A}\rightarrow [0,\infty)$ defined by 
$\| x \|_{a,\lambda}^2 := \sup\{\lambda \varphi (x^* a x ) + (1-\lambda ) | \varphi (ax)|^2  :  \varphi \in \mathcal{S}_a (\mathcal{A}) \}$ for all $ \lambda \in [0,1] $.
These family of norms generalize the well-known norms $\|\cdot\|_a$ and $v_a(\cdot)$ on $\mathcal{A}$. The notions of norm-parallelism and  Birkhoff-James orthogonality with respect to $\|\cdot\|_{a,\lambda}$, which is called, $(a,\lambda)$-norm parallelism and $(a,\lambda)$-Birkhoff-James orthogonality respectively, are introduced and investigated. We characterize $(a,\lambda)$-norm parallelism and $(a,\lambda)$-Birkhoff-James orthogonality in terms of the elements of $\mathcal{S}_a(\mathcal{A})$. As special cases, characterization of Birkhoff-James orthogonality and norm parallelism with respect to $\|\cdot\|_a$ and $v_a(\cdot)$ are obtained. In particular, the connection between these new concepts are described. The
results in this paper are extensions and generalizations of known
results which appeared in \cite{EAZ,Moslehian,ZM1,ZM2,Z2,ZW,AZ}.

%%%%%%%%%%%%%%%%%%%%%%%%%%%%%%%%%%%%%%%%%%%
%\section{$(a,\lambda)$-norm parallelism in unital $C^*$-algebras}
\section{$(a,\lambda)$-semi-norms on unital $C^*$-algebras}
We start this section by defining a family of semi-norms as a generalization of semi-norms $\|\cdot\|_a$ and $v_a(\cdot)$ on unital $C^*$-algebras.
\begin{definition}
Let $ \mathcal{A} $ be a unital $ C^* $-algebra, $ a \in \mathcal{A} $ be a positive element and $ \lambda \in [0,1] $. For every $ x \in \mathcal{A} $, define
\begin{align*}
\| x \|_{a,\lambda }:= \sup \{ \sqrt{\lambda \varphi (x^* a x ) + (1-\lambda ) | \varphi (ax)|^2}  :  \varphi \in \mathcal{S}_a (\mathcal{A}) \}.
\end{align*}
\end{definition}
Note that for every $ x \in \mathcal{A} $, we have $ \| x \|_{a,1}= \| x \|_{a} $ and $ \| x \|_{a,0} = v_a (x) $. Also,  moving $\lambda$ along the interval $[0, 1]$, a wide of seminorms are presented.

In this section, we collect some main properties on $(a,\lambda)$-semi-norm on unital $C^*$-algebras and its relation with special semi-norms $\|\cdot\|_a$ and $v_a(\cdot)$.

\begin{proposition}\label{007}
The functional $ \| \cdot \|_{a,\lambda } $ is a semi-norm on $ \mathcal{A} $ and the following inequalities always hold for every $ x \in \mathcal{A} $:
\begin{align}\label{001}
v_a (x) \leq \| x \|_{a,\lambda } \leq \| x \|_{a} .
\end{align}
\end{proposition}
\begin{proof}
Let $ x \in \mathcal{A} $ . Obviously, $ \| x \|_{a,\lambda } \geq 0 $. For each $ \alpha \in \mathbb{C} $, we have
\begin{align*}
\| \alpha x \|_{a,\lambda } &= \sup \{ \sqrt{\lambda \varphi (\overline{\alpha} x^* a \alpha x ) + (1 - \lambda ) | \varphi (a \alpha x )|^2 } \, :\, \varphi \in \mathcal{S}_a (\mathcal{A})  \} \\
&= \sup \{ \sqrt{| \alpha |^2 \lambda \varphi ( x^* a x ) +|\alpha |^2 (1 - \lambda ) | \varphi (a x )|^2 } \, :\, \varphi \in \mathcal{S}_a (\mathcal{A})  \} \\
&= \sup \{ |\alpha | \sqrt{\lambda \varphi (x^* a x ) + (1 - \lambda ) | \varphi (a x )|^2 } \, :\, \varphi \in \mathcal{S}_a (\mathcal{A}) \} \\
&= | \alpha | \| x \|_{a , \lambda} .
\end{align*}
Now, let $ \varphi \in \mathcal{S}_a (\mathcal{A}) $, $ x,y \in \mathcal{A} $ and  $ \lambda \in [0,1] $. Then
\begin{align*}
&\lambda | \varphi ((x+ y )^* a (x+y) ) +(1-\lambda ) | \varphi (a(x+y))|^2 \\
= \,  &\lambda\big( \varphi (x^* a x) + 2 {\rm Re} \varphi (x^* a y ) +\varphi (y^* a y)\big) \\
&+(1- \lambda )\big( | \varphi(ax)|^2 + 2 {\rm Re}  \overline{\varphi (ax)} \varphi (ay) + | \varphi (ay)|^2\big) \\
\leq \, &\lambda \big(\varphi (x^* a x) + 2 | \varphi (x^* a y )| +  \varphi (y^* a y)\big) \\
&+(1- \lambda ) \big(| \varphi(ax)|^2 + 2| \varphi (ax)| | \varphi (ay)| + | \varphi (ay)|^2 \big) \\
\leq \, &\lambda \big(\varphi (x^* a x) + 2 \varphi (x^* a x )^{\frac{1}{2}} \varphi (y^* a y )^{\frac{1}{2}} +  \varphi (y^* a y)\big) \\
&+(1- \lambda )\big( | \varphi (ax)|^2 + 2| \varphi (ax)| | \varphi (ay)| + | \varphi (ay)|^2\big) \\
& \text{(by Cauchy-Schwartz inequality)}\\
= \,&(\lambda \varphi (x^* a x) + (1-\lambda) |\varphi (ax)|^2) + (\lambda \varphi (y^* a y)+ (1-\lambda) |\varphi (ay)|^2)\\
&+  2\big( \lambda  \varphi (x^* a x )^{\frac{1}{2}} \varphi (y^* a y )^{\frac{1}{2}} + (1-\lambda) |\varphi (ax)||\varphi (ay)| \big)\\
\leq \, & (\lambda \varphi (x^* a x) + (1-\lambda) |\varphi (ax)|^2 )+( \lambda \varphi (y^* a y)+ (1-\lambda) |\varphi(ay)|^2)  \\
&+ 
2 \sqrt{\lambda \varphi (x^* a x) + (1-\lambda) |\varphi (ax)|^2 } \sqrt{\lambda \varphi (y^* a y) + (1-\lambda) |\varphi (ay)|^2 })\\
&\text{(by Cauchy-Bunyakovsky-Schwartz inequality)} \\
\leq \,  &\| x  \|_{a , \lambda }^{2} + \| y \|_{a , \lambda }^{2} +  2 \| x \|_{a , \lambda } \| y \|_{a , \lambda } 
= (\| x  \|_{a , \lambda } + \|  y \|_{a , \lambda } )^2 .
\end{align*}
Taking a supremum on all $ \varphi \in \mathcal{S}_a (\mathcal{A}) $, we obtain
\begin{align*}
\| x + y \|_{a , \lambda }^{2} \leq  (\| x  \|_{a , \lambda } + \|  y \|_{a , \lambda } )^2,
\end{align*}
which follows the triangle inequality 
$ \| x + y \|_{a , \lambda } \leq \| x \|_{a , \lambda } + \| y \|_{a , \lambda } $.

For the proof of \eqref{001}, let $ \varphi \in \mathcal{S}_a (\mathcal{A}) $. By the Cauchy-Schwartz inequality, we have
\begin{align*}
|\varphi (ax)|^2 = \lambda |\varphi (ax)|^2  + (1-\lambda ) |\varphi (ax)|^2
 \leq \lambda \varphi (a) \varphi (x^* a x) + (1-\lambda ) |\varphi (ax)|^2.
\end{align*}
Hence
$$ |\varphi (ax)| \leq \sqrt{\lambda  \varphi (x^* a x) + (1-\lambda ) |\varphi (ax)|^2} . $$
Therefore
$$ v_a (x) = \sup_{\varphi \in \mathcal{S}_a (\mathcal{A}) } | \varphi (ax)|  \leq \| x \|_{a,\lambda }. $$
Moreovere, 
\begin{align*}
\lambda \varphi (x^* a x) + (1-\lambda ) | \varphi (ax)|^2 
\leq  \lambda \varphi (x^* a x) + (1-\lambda) \varphi (a) \varphi (x^* a x) = \varphi (x^* a x).
\end{align*}
It follows that 
\[\|x\|_{a,\lambda}=\sup_{\varphi \in \mathcal{S}_a (\mathcal{A}) } \sqrt{\lambda  \varphi (x^* a x) + (1-\lambda ) |\varphi (ax)|^2}\leq\sup_{\varphi \in \mathcal{S}_a (\mathcal{A})}\sqrt{\varphi(x^*ax)}=\|x\|_a.\]
\end{proof}
In \cite[Corollary 4.5]{AB}, was shown that if $ x \in \mathcal{A} $ is an $ a $-selfadjoin element, then $ v_a (x)=\| x\|_a $. On the other hand, by \eqref{001}, we have $ v_a (x) \leq \| x\|_{a ,\lambda} \leq \|x\|_a $. Thus
$$  v_a (x)=\| x\|_a \leq \| x\|_{a ,\lambda} \leq \|x\|_a, $$
and so
\begin{align}\label{008}
 v_a (x)=\| x\|_a=\| x\|_{a ,\lambda}
\end{align}
for all $ a $-selfadjoint elements $ x \in \mathcal{A} $.
Also, if $x\in\mathcal{A}_a$, then $x=x_1+ix_2$, where $x_1=\dfrac{x+x^{\sharp}}{2}$ and $x_2=\dfrac{x-x^{\sharp}}{2i}$, and therefore
\begin{align*}
\|x\|_{a,\lambda}=\| x_1+ i x_2\|_{a,\lambda}\leq \| x_1\|_{a,\lambda} +\|  x_2\|_{a,\lambda}
= v_a( x_1) + v_a( x_2) \leq v_a (x)+v_a(x) =2v_a (x).
\end{align*}
Thus \[ \dfrac{1}{2} \|x\|_{a,\lambda}  \leq v_a(x) \leq \|x\|_{a,\lambda}.\]
%%%%%%%%%%%%%%%%%%%%%%%%%%%%%%%%%%%%%%
Let $ \mathcal{H}$ be a complex Hilbert space with the inner product $ \langle \cdot , \cdot \rangle $ and let $ \mathcal{B}(\mathcal{H}) $ be the $ C^* $-algebra of all bounded linear operators on $ \mathcal{H}$. Assume that $A\in\mathcal{B}(\mathcal{H})^+$, which induces a positive semi-definite sesquilinear form
$\langle\cdot,\cdot\rangle_A:\mathcal{H}\times\mathcal{H}\rightarrow\mathbb{C}$ defined by $\langle x,y\rangle_A=\langle Ax,y\rangle$.
The semi-norm $\|\cdot\|_A$ induced by $\langle\cdot,\cdot\rangle_A$ is defined by $\|x\|_A=\sqrt{\langle Ax,x\rangle}$ for every $x\in \mathcal{H}$; cf. \cite{AM}.
 The set of all $A$-bounded operators on $\mathcal{H}$ is defined by
$$ \mathcal{B}_{A^{\frac{1}{2}}} ( \mathcal{H}) := \{ T \in \mathcal{B}(\mathcal{H}) \, : \, \exists \, c > 0: \, \| T x \|_{A} \leq c \| x \|_{A}, \quad \forall x \in \mathcal{H} \}.$$
In fact, $\mathcal{B}_{A^{\frac{1}{2}}} ( \mathcal{H})$ is a unital subalgebra of $\mathcal{B}(\mathcal{H}) $ which is equipped with the semi-norm
$$ \gamma_A (T) := \sup_{\|x\|_A=1} \sqrt{\langle ATx , Tx \rangle }\quad(T \in \mathcal{B}_{A^{\frac{1}{2}}} ( \mathcal{H})). $$
It was proved in \cite{Mabrouk,AB} that if $\pi :\mathcal{A}\rightarrow\mathcal{B}(\mathcal{H})$ is a unital faithful $*$-representation of $\mathcal{A}$, then
\begin{align}\label{01009}
\| x \|_a = \gamma_{\pi(a)} ( \pi(x)),\quad v_a(x)=v_{\pi(a)}(\pi(x)) \quad(x\in \mathcal{A}).
\end{align}
Let $ \lambda \in [ 0 ,1 ] $. For every $ T \in \mathcal{B}_{A^{\frac{1}{2}} } (\mathcal{H}) $ the semi-norm
\begin{align}\label{1300}
\gamma_{A , \lambda } (T) =\sup \{ \sqrt{\lambda \| T \|_{A}^2 + (1-\lambda )| \langle Tx , x \rangle_{A}|^2 }:\,\, x \in \mathcal{H}, \, \| x\|_A =1 \}, 
\end{align}
is a generalization of $\gamma_A(\cdot)$ which is introduced and studied in \cite{EAZ}.

A positive linear functional $f$ on $\mathcal{A}$ is said to be pure if for every positive functional $ g $ on $\mathcal{A}$ satisfying $  g(xx^*
) \leq f (xx^*) $ for all $ x\in\mathcal{A} $, there is a scalar $ 0\leq \lambda \leq 1 $ such that $ g=\lambda f $. The set of pure states on $ \mathcal{A} $ is denoted by $ \mathcal{P}(\mathcal{A}) $.

In the next theorem, we prove that there is a relationship similar to \eqref{01009} between  $ \| \cdot \|_{a , \lambda} $ and $ \gamma_{A , \lambda}(\cdot) $. The proof method employed in this theorem follows the techniques used in Theorem 3.5 of \cite{AB}. So, we need the following known lemmas.
\begin{lemma}\cite[corollary 5.1.10]{Murphy}\label{1500}
Let $ \mathcal{A} $ be a unital $ C^* $-algebra. Then $ \mathcal{S}(\mathcal{A}) $
is the weak$^* $-closed convex hull of the space of pure states on $\mathcal{A} $.
\end{lemma}
\begin{lemma}\cite[lemma A.41]{RIW}\label{1400}
Suppose that $ \mathcal{F} $ is a subset of $ \mathcal{S}(\mathcal{A}) $ such that the following statement holds true: an element $ x \in \mathcal{A} $ is positive if and only if $ x $ is selfadjoin and $ \rho(x) \geq 0 $ for all $ \rho \in \mathcal{F} $. Then $ \mathcal{P}(\mathcal{A}) $ is contained in the $ weak^* $-closure of $\mathcal{ F} $. 
\end{lemma}

\begin{theorem}\label{160}
Let $ \mathcal{A} $ be a unital $ C^* $-algebra, $ a\in \mathcal{A} $ be positive and let $ \pi : \mathcal{A} \rightarrow \mathcal{B}(\mathcal{H}) $ be a unital faithful $ * $-representation of $ \mathcal{A}$. Then 
$$ \| x\|_{a,\lambda } = \gamma_{\pi(a) , \lambda }  (\pi (x)) $$
for all $ x \in \mathcal{A}^a $ and all $ \lambda \in [0,1] $.
\end{theorem}
\begin{proof}
Suppose that $ h\in \mathcal{H} $ is a $\pi(a)$-unit vector; i.e., $ \langle \pi (a) h , h \rangle = 1 $. The functional $ \varphi_h : \mathcal{A} \rightarrow \mathbb{C} $ defined by $ \varphi_h (z) = \langle \pi (z) h , h \rangle  $ belongs to $ \mathcal{S}_{a} (\mathcal{A}) $. So for each $ \lambda \in [0,1]$, we have
\begin{align*}
&\sqrt{\lambda \| \pi (x ) h \|_{\pi (a)}^{2} + (1-\lambda ) | \langle \pi (x)h , h \rangle_{\pi (a)} |^2 }  \\
= & \sqrt{\lambda  \langle \pi (a) \pi (x)h ,\pi (x) h \rangle +  (1-\lambda) |  \langle \pi (a) \pi (x)h , h \rangle |^2 } \\ 
= &\sqrt{\lambda  \langle \pi (x^* a x ) h , h \rangle +  (1-\lambda) | \langle \pi (ax) h , h \rangle |^2 } \\
=&\sqrt{\lambda \varphi_h (x^* ax ) + (1-\lambda) | \varphi_h (ax ) |^2 } \\
\leq &\sup_{\varphi \in \mathcal{S}_a (\mathcal{A}) }\sqrt{\lambda \varphi (x^* ax ) + (1-\lambda) | \varphi (ax ) |^2 } = \| x \|_{a , \lambda }
\end{align*}
Taking supremum on all $ h \in \mathcal{H} $ with $ \| h \|_{\pi (a)} = 1 $, we obtain
\begin{align*}
 \gamma _{\pi (a) , \lambda } ( \pi (x)) = \sup_{\| h \|_{\pi (a)} = 1} \sqrt{\lambda \| \pi (x ) h \|_{\pi (a)}^{2} + (1-\lambda ) | \langle \pi (x)h , h \rangle_{\pi (a)} |^2}  \leq \| x \|_{a , \lambda}<+\infty,
\end{align*} 
since $ \| x \|_{a , \lambda } \leq \| x \|_{a }<+\infty $.

Now, let $ h \in \mathcal{H} $ be arbitrary. Then $\langle \pi (a) h , h \rangle\neq0$ and $ \dfrac{h}{\sqrt{\langle \pi (a) h , h \rangle }}$ is a $\pi(a)$-unit vector of $\mathcal{H}$.  So, we have
 \begin{align}\label{1600}
 \lambda \dfrac{1}{\langle \pi (a) h , h \rangle } \langle \pi ( x^* a x ) h , h \rangle + \dfrac{1-\lambda}{\langle \pi (a) h , h \rangle^2} | \langle \pi (ax ) h,h \rangle |^2 \leq \gamma^2_{\pi ( a),\lambda } (\pi  (x)) .
\end{align}
Let $ \rho_h $ be the vector state defined by $ \rho_h (T):= \langle T h ,h \rangle $ for all $  T \in \mathcal{B}(\mathcal{H}) $. Set $ \mathcal{F}:=\{ \rho_h \circ \pi \, : \, h \in \mathcal{H} \} $. It is trivial that $\rho_h\circ\pi$ is linear for all $h\in\mathcal{H}$ and 
$$  \rho_h \circ \pi (1_\mathcal{A})=\rho_h(I_\mathcal{H})=\|h\|^2=\| \rho_h \circ \pi \| =1. $$
Hence $\rho_h\circ\pi\in\mathcal{S}(\mathcal{A}) $ for all $h\in\mathcal{H}$, and so $\mathcal{F} \subset \mathcal{S}(\mathcal{A}) $.
Thus, the conditions of Lemma \ref{1400}  hold. Therefore $ \mathcal{S}(\mathcal{A}) = \overline{\mathcal{F}}^{w^*} $.
Assume that  $ f $ is a pure state on $ \mathcal{A} $. It follows from Lemma \ref{1500} that, there exists a net $ \{ h_\alpha \} _{\alpha \in \Lambda } \subset \mathcal{H} $ such that $ \{ \rho_{h_\alpha } \circ \pi \}_{\alpha \in \Lambda } $ weak$^* $ converges to $ f $. By replacing net $ \{ h_\alpha\}_{\alpha \in \Lambda } $ in \eqref{1600} we get
\begin{align*}
\dfrac{\lambda \langle \pi (x^* a x) h_\alpha , h_\alpha \rangle}{\langle\pi (a) h_\alpha , h_\alpha \rangle } + \dfrac{(1-\lambda)}{\langle\pi (a) h_\alpha , h_\alpha \rangle^2}|\langle\pi (ax) h_\alpha , h_\alpha \rangle |^2 \leq \gamma_{\pi(a),\lambda}^{2} (\pi(x)).
\end{align*}
Hence
$$ \dfrac{\lambda \rho_{h_\alpha} \circ \pi (x^* ax)}{\rho_{h_\alpha} \circ \pi ( a)} + \dfrac{(1-\lambda)}{\rho_{h_\alpha} \circ \pi (a)} | \rho_{h_\alpha} \circ \pi ( ax) |^2 \leq \gamma_{\pi(a),\lambda}^{2} (\pi(x)). $$
Taking limit from the above inequality, infer that
\begin{align}\label{Eqf}
\dfrac{\lambda f (x^* ax)}{f ( a)} + \dfrac{(1-\lambda)}{f (a)^{2}} | f( ax) |^2 \leq \gamma_{\pi(a),\lambda}^{2} (\pi(x))\quad(\forall f \in \mathcal{P}(\mathcal{A})).
\end{align}
Moreover, applying the same technique used in the proof of Theorem 3.5 of \cite{AB} and using Lemma \ref{1400} ensure that the inequality \eqref{Eqf} holds true for all $f\in\mathcal{S}(\mathcal{A})$.  So, we conclude that
$$\lambda g (x^*  ax) + (1 - \lambda ) | g(ax)|^2 \leq  \gamma_{\pi(a),\lambda}^{2} (\pi(x))\quad(\forall g\in  \mathcal{S}_a (\mathcal{A})). $$
Taking supremum on all $  g \in \mathcal{S}_a (\mathcal{A}) $, follows that
$$ \| x \|_{a , \lambda } \, \leq  \gamma_{\pi(a),\lambda} (\pi(x)), $$
and the proof is completed.
\end{proof}
%%%%%%%%%%%%%%%%%%%%%%%%%%%%%%%%%%%%%%%%%%%%%%%%%%%%
\begin{corollary}\label{009}
Let $ \mathcal{A} $ be a unital $ C^* $-algebra, $ a \in \mathcal{A} $ be a positive element and let $ x \in \mathcal{A}_a$. Then
\begin{align}\label{CS}
\| x\|_{a,\lambda } = \| x^\sharp\|_{a,\lambda} \quad (\forall \lambda \in [0,1]).
\end{align}
\end{corollary}
\begin{proof}
Let $ x \in \mathcal{A}_a$ and $ h\in \mathcal{H} $ be such that $ \|h\|_{\pi(a)}=1 $.  By \eqref{s0009} and \eqref{01009}, we have
$$ \gamma_{\pi(a)} (\pi (x))=\|x\|_a^2 =\|x^{\sharp}\|_a^2 = \gamma_{\pi(a)} (\pi (x^{\sharp})). $$
So Theorem \ref{160} yields that
\begin{align*}
\|x\|_{a,\lambda}^2 =\gamma_{\pi(a),\lambda}^2 (\pi (x))&= \sup_{\|h\|_{\pi(a)}=1} ( \lambda \gamma_{\pi(a)} (\pi(x)) + (1-\lambda) | \langle \pi(x) h,h\rangle_{\pi(a)}|^2 ) \\
&= \sup_{\|h\|_{\pi(a)}=1} ( \lambda \gamma_{\pi(a)} (\pi(x^{\sharp})) + (1-\lambda) | \langle \pi(ax) h,h\rangle |^2 ) \\
&= \sup_{\|h\|_{\pi(a)}=1} ( \lambda \gamma_{\pi(a)} (\pi(x^{\sharp})) + (1-\lambda) | \langle \pi(x^* a) h,h\rangle |^2 )\\
&= \sup_{\|h\|_{\pi(a)}=1} ( \lambda \gamma_{\pi(a)} (\pi(x^{\sharp})) + (1-\lambda) | \langle \pi(ax^{\sharp}) h,h\rangle |^2 ) \\
&= \sup_{\|h\|_{\pi(a)}=1} ( \lambda \gamma_{\pi(a)} (\pi(x^{\sharp})) + (1-\lambda) | \langle \pi(x^{\sharp}) h,h\rangle_{\pi(a)} |^2 )\\
&=\gamma_{\pi(a),\lambda}^2 (\pi (x^{\sharp})) =\|x^{\sharp}\|_{a,\lambda}^2.
\end{align*}
\end{proof}
%%%%%%%%%%%%%%%%%%%%%%%%%%%%%%%%%%%%%%%%%%%%%%%%%%%%%%%%%
Let $ x \in \mathcal{A}_a $ and $ x^{\sharp} $ be an $ a $-adjoint of $x$. Since $ x x^{\sharp} $ is $ a $-selfadjoint, from \eqref{s0009} and \eqref{008} we get
\begin{align}\label{0077}
 \|x x^{\sharp}\|_{a,\lambda}= \|x x^{\sharp}\|_{a}=\|x \|_{a}^2 .
\end{align}
So $  \|x \|_{a,\lambda}^2 \leq  \|x x^{\sharp}\|_{a,\lambda} $, by \eqref{001}. 
But the equality $\|xx^{\sharp}\|_{a,\lambda} = \|x\|^2_{a,\lambda}$ is not established
in general. The following example illustrate this fact.
\begin{example}\label{1100}
 Let ${\rm Tr}$ be the usual trace functional on $ C^* $-algebra of all $ 2 \times 2 $ complex matrices $ \mathbb{M}_2 (\mathbb{C}) $ with identity matrix $I_2$ as unit. According to the Example 2.2 of \cite{AB}, for any positive matrix $  h \in \mathbb{M}_2 (\mathbb{C}) $, let $\varphi_h$ be the positive linear functional given by 
 \[\varphi_h (x)={\rm Tr}(hx) \quad (x\in \mathbb{M}_2 (\mathbb{C})).\]
Then for a positive matrix $ a\in \mathbb{M}_2 (\mathbb{C}) $,  we have
\begin{align*}
\mathcal{S}_{a} (\mathbb{M}_{2}(\mathbb{C})) = \{ \varphi_h \, : \, h \in  \mathbb{M}_2 (\mathbb{C})^{+} \quad {\rm and} \quad {\rm Tr}(ha) =1\}.
\end{align*}
Let $ a=\left[\begin{array}{cccc}
2&0  \\
0&1
\end{array}\right]  $. Then with some simple matrix computations, we conclude that
\begin{align*}
\mathcal{S}_{a} (\mathbb{M}_{2}(\mathbb{C})) =\{ \varphi_h \, : \, h\in\mathcal{L}_a\},
\end{align*}
where
{\small\begin{align*}
\mathcal{L}_a:=\{h
=\left[\begin{array}{cccc}
h_{11}&h_{12}  \\
\overline{h}_{12}& h_{22}
\end{array}\right]
\in  \mathbb{M}_2 (\mathbb{C})^{+}\, : \, h_{12}  \in \mathbb{C}, \, h_{11}, h_{22} \geq 0 \,\,\, {\rm and} \,\,\, 2h_{11} +h_{22} =1\}.
\end{align*}}
Let $
x=\left[\begin{array}{cccc}
0&1  \\
0&0
\end{array}\right]
$. By  \eqref{0077}, we have
\begin{align*}
&\|x x^{\sharp}\|_{a,\lambda} = \| x \|_{a}^{2} =
\sup_{\varphi_{h} \in \mathcal{S}_a (\mathbb{M}_2 (\mathbb{C}))} \varphi_{h} (x^* a x) \\
&= \sup_{h\in\mathcal{L}_a} {\rm Tr}\big(
\left[\begin{array}{cccc}
0& 0  \\
0& 2h_{22}
\end{array}\right]
\big)
 =\sup_{2h_{11}+h_{22}=1, h_{11}, h_{22} \geq 0} 2h_{22} = 2.
\end{align*}
Note that $h
=\left[\begin{array}{cccc}
h_{11}&h_{12}  \\
\overline{h}_{12}& h_{22}
\end{array}\right]\in\mathbb{M}_2 (\mathbb{C})^+$ if and only if $h_{11}\geq0$, $h_{22}\geq0$ and there is $k\in\mathbb{C}$ with $| k|\leq1$ such that $h_{12}=\sqrt{h_{11}h_{22}}k$. Hence for any $ \lambda \in [0,1) $, we have
\begin{align*}
\|x\|_{a,\lambda}^2 &=\sup_{\varphi_{h} \in \mathcal{S}_a (\mathbb{M}_2 (\mathbb{C}))} (\lambda  \varphi_{h} (x^* a x)  +(1-\lambda)|\varphi_{h} (a x)|^2 ) \\
&=\sup_{h\in\mathcal{L}_a} (2\lambda h_{22} + 4(1-\lambda) | h_{12}|^2 ) \\
&=\sup_{2h_{11}+h_{22}=1, h_{11}, h_{22} \geq 0,|k|\leq1} (2\lambda h_{22} + 4(1-\lambda) (h_{11} h_{22}) |k|^2 ) \\
 & \leq \, 2\sup_{2h_{11}+h_{22}=1, h_{11}, h_{22} \geq 0} (\lambda h_{22} + 2(1-\lambda) h_{11} h_{22} ) \\
%&= \sup_{h_{22}\in [0,1]} (\lambda 2h_{22} + 4(1-\lambda) (\dfrac{1-h_{22}}{2}) h_{22}) \\
&= 2\sup_{h_{22}\in [0,1]} ( h_{22} - (1-\lambda) h^2_{22}).
\end{align*}
Note that  $\sup_{h_{22}\in[0,1]} (h_{22} - (1-\lambda) h^2_{22} )=\frac{1}{4(1-\lambda)}$ for all $\lambda\in [0,1)$. Hence, if we assume that $\lambda\in [0,\frac{3}{4})$, then we have $\sup_{h_{22}\in[0,1]} (h_{22} - (1-\lambda) h^2_{22} )<1$, and so
\begin{align*}
\|x\|_{a,\lambda}^2 \leq2\sup_{h_{22}\in [0,1]}( h_{22} - (1-\lambda) h_{22}^2 )<2=\| x\|^2_a=\|xx^{\sharp}\|^2_{a,\lambda}.
\end{align*}
\end{example}
%%%%%%%%%%%%%%%%%%%%%%%%%%%%%%%%%%%%%%%%%%%%%%%%%%%%555
%%%%%%%%%%%%%%%%%%%%%%%%%%%%%%%%%%%%%%%%%%%%%%%%%%%%%%%%%
%%%%%%%%%%%%%%%%%%%%%%%%%%%%%%%%%%%%%%%%%%
\section{Characterization of  $ (a,\lambda) $-norm parallelism in unital $C^*$-algebra}
Let $\mathcal{A}$ be a unital $C^*$-algebra and let $a\in\mathcal{A}$ be a nonzero positive element. It was proved in \cite[Theorem 3.9]{AB} that $\mathcal{A}_a\subset\mathcal{A}^a$.  Now, if we assume that $a$ is a positive and invertible element of $\mathcal{A}$, then for every $x\in\mathcal{A}$ the equation $ay=x^*a$ has the unique solution $x^{\sharp}=a^{-1}x^*a$, and so every $x\in\mathcal{A}$ is $a$-adjointable. Therefore $\mathcal{A}^a=\mathcal{A}$.
From now on we assume that $\mathcal{A}$ is a unital $C^*$-algebra and $a\in\mathcal{A}$ is a positive and invertible element.

First of all, note that if $ a\in \mathcal{A} $ is positive and invertible, then $ \| \cdot \|_{a,\lambda} $ is a norm on $\mathcal{A}$ for all $ \lambda \in [0,1] $. Indeed, if  $x \in \mathcal{A}$ and $ \| x\|_{a,\lambda }=0 $, then $ \lambda \varphi (x^* a x ) + (1-\lambda ) | \varphi (ax)|^2 =0 $ for all $ \varphi \in \mathcal{S}_a (\mathcal{A}) $ and all $ \lambda \in [0,1] $. Then $ \varphi (x^* a x ) =\varphi (ax)=0 $ for all $ \varphi \in \mathcal{S}_a (\mathcal{A}) $, which implies that $ ax=0 $, 
and so $x=0$, since $a$ is invertible. 

We start this section by defining a new kinds of norm-parallelism in unital $C^*$-algebras.
\begin{definition}
Let $ \mathcal{A} $ be a unital $ C^*$-algebra, $a\in\mathcal{A}$ be a positive and invertible element and $\lambda\in[0,1]$. An element $ x \in \mathcal{A} $ is called $ (a,\lambda )$-norm parallel to $ y\in \mathcal{A} $, denoted by $ x ||_{a,\lambda} y $, if there exists $ \mu \in \mathbb{T}:=\{z \in \mathbb{C}\, : \, | z |=1\} $ such that
$$ \| x + \mu y \|_{a, \lambda } = \| x\|_{a ,\lambda}+\| y\|_{a ,\lambda}. $$
\end{definition}
Note that norm parallelism and numerical radius parallelism with respect $\|\cdot\|_a$ and $v_a(\cdot)$ is just $(a,1)$-norm parallelism and $(a,0)$-norm parallelism, resectively which can be defined naturally as follows:

$\bullet$ An element $x\in\mathcal{A}$ is called  $a$-norm parallel to
another element $y\in\mathcal{A}$, denoted by $x\parallel_a y$,  if there is $\mu\in\mathbb{T}$ such that $\|x+\mu y\|_a=\|x\|_a+\|y\|_a$.

$\bullet$ An element $x\in\mathcal{A}$ is called $a$-numerical radius parallel to
another element $y\in\mathcal{A}$, denoted by $x\parallel_{v_a}y$, if there is $\mu\in\mathbb{T}$ such that $v_a(x+\mu y)= v_a(x)+v_a(y)$.

The following proposition discusses the main properties of $(a,\lambda)$-parallelism.
The proof follows by properties of $(a,\lambda)$-norm and by using the idea of the proof of Proposition 3.2
in \cite{ZM1}; and so the details are omitted.
\begin{proposition}\label{16}
Let $ \mathcal{A} $ be a unital $ C^*$-algebra, $a\in\mathcal{A}$ be a positive and invertible element. If $ x,y \in \mathcal{A}$, then the following statements are equivalent:
\begin{itemize}
\item[1)] $ x ||_{a,\lambda} y $;
\item[2)] $ x^\sharp ||_{a,\lambda} y^\sharp $;
\item[3)] $ \alpha x ||_{a,\lambda} \beta y $ for any $ \alpha ,\beta \in \mathbb{R}\setminus \{0\} $;  i.e., $ (a,\lambda )$-norm parallelism is $\mathbb{R}$-homogenous.
\end{itemize}
\end{proposition}
%%%%%%%%%%%%%%%%%%%%%%%%%%%%%%%%%%%%%%%%%
%%%%%%%%%%%%%%
In the following result we give a characterization of $(a,\lambda)$-parallelism of elements of $C^*$-algebra $\mathcal{A}$ in terms of elements of $\mathcal{S}_a({\mathcal{A}})$.
%%%%%%%%%%%%%%%%%%
\begin{theorem}\label{002}
Let $ \mathcal{A} $ be a unital $ C^* $-algebra, $ a \in \mathcal{A} $ be a positive and invertible element and $ x,y \in \mathcal{A} $. Then the following statements are equivalent:
\begin{itemize}
\item[1)] $ x \|_{a,\lambda} y $; i.e., $ \| x+\mu y \|_{a , \lambda } = \| x \|_{a,\lambda } + \| y \|_{a , \lambda } $ for some $ \mu \in \mathbb{T} $;
\item[2)] There exist $ \varphi \in \mathcal{S}_a (\mathcal{A}) $ and $ \mu \in \mathbb{T} $ such that
\begin{align}
 \lambda \varphi (x^* ay) +(1-\lambda ) \overline{\varphi(ax)} \varphi (ay)=\overline{\mu} \| x \|_{a , \lambda } \,  \| y \|_{a , \lambda }.
\end{align}
\item[3)] $ \| x \|_{a,\lambda}  \|y \|_{a,\lambda} \in \lambda V_a (x^\sharp y)
+ (1-\lambda) V_a (x^\sharp) V_a (y).$
\end{itemize}
\end{theorem}
\begin{proof}
First note that the implication $2)\Leftrightarrow 3)$ are clear.
\item[$ 1) \Rightarrow 2 $)]
Assume that $ \| x+\mu y\|_{a,\lambda}=\| x\|_{a,\lambda}+\| y\|_{a,\lambda}$ for some $ \mu \in \mathbb{T} $.
By \cite[Proposition 2.3]{AB},  $ \mathcal{S}_a (\mathcal{A}) $ is $ w^* $-compact, since $ a $ is invertible. So there exists $ \varphi \in \mathcal{S}_a (\mathcal{A}) $ such that
{\small\begin{align*}
&(\| x\|_{a,\lambda}+\| y\|_{a,\lambda})^2 =\| x+\mu y\|_{a,\lambda}^2\\
&= \lambda \varphi ((x+\mu y)^* a (x+\mu y))+(1-\lambda) |\varphi(a(x+\mu y))|^2 \\
 &= \, \lambda \big(\varphi (x^* a x)+|\mu |^2 \varphi (y^* a y)+2 \mathrm{Re}(\mu \varphi (x^* ay))\big)\\
 &+(1-\lambda) \big(|\varphi(ax)|^2 +|\mu |^2 |\varphi(ay)|^2 +2 \mathrm{Re} (\mu \overline{\varphi(ax)} \varphi(ay))\big)\\
= & \, (\lambda \varphi (x^* a x) + (1-\lambda) |\varphi(ax)|^2 ) + (\lambda \varphi (y^* a y) +(1-\lambda) |\varphi(ay)|^2 ) \\
  &+ 2\mathrm{Re} \big(\lambda\mu \varphi (x^* ay)
 + (1-\lambda)\mu \overline{\varphi(ax)} \varphi(ay)\big) \\
 \leq &  \, \| x\|_{a,\lambda}^2 +\| y\|_{a,\lambda}^2 + 2\mathrm{Re} \big(\lambda\mu \varphi (x^* ay)
 + (1-\lambda)\mu \overline{\varphi(ax)} \varphi(ay)\big) \\
 &\leq \,   \| x\|_{a,\lambda}^2 +\| y\|_{a,\lambda}^2 + 2 \lambda | \varphi (x^* ay)|+ 2(1-\lambda)  |\varphi(ax)|\,|\varphi(ay)| \\
 &\leq \,  \| x\|_{a,\lambda}^2 +\| y\|_{a,\lambda}^2 + 2\lambda \varphi(x^* ax)^{\frac{1}{2}}\varphi(y^* ay)^{\frac{1}{2}} +
 2(1-\lambda) |\varphi(ax)| |\varphi(ay)|\\
 &\quad \text{(by the Cauchy-Schwartz inequality)} \\
 \leq &\, \| x\|_{a,\lambda}^2 +\| y\|_{a,\lambda}^2 + 2
\sqrt{\lambda \varphi (x^* a x)+(1-\lambda) |\varphi(ax)|^2 }
\sqrt{\lambda \varphi (y^* a y)+(1-\lambda) |\varphi(ay)|^2 }\\
&\quad \text{(by the Cauchy-Bunyakovsky-Schwartz inequality)}\\
\leq &\,  \| x\|_{a,\lambda}^2 +\| y\|_{a,\lambda}^2 + 2 \| x\|_{a,\lambda} \| y\|_{a,\lambda} =(\| x\|_{a,\lambda}+\| y\|_{a,\lambda})^2 .
\end{align*}}
Thus there exists $ \varphi \in \mathcal{S}_a (\mathcal{A}) $ such that
\begin{align}\label{32}
\mathrm{Re} (\lambda \mu \varphi(x^* ay) + (1-\lambda)\mu \overline{\varphi(ax)}\varphi(ay) )=
\| x\|_{a,\lambda}\| y\|_{a,\lambda}.
\end{align}
Moreover, for every $ \lambda\in [0,1] $, we have
{\small\begin{align*}
\mathrm{Re}^2 \big(\lambda \mu \varphi(x^* ay) &+ (1-\lambda)\mu \overline{\varphi(ax)}\varphi(ay) \big)
+ \mathrm{Im}^2 \big(\lambda \mu \varphi(x^* ay) + (1-\lambda)\mu \overline{\varphi(ax)}\varphi(ay) \big)\\
&= |\lambda \mu \varphi(x^* ay) + (1-\lambda)\mu \overline{\varphi(ax)}\varphi(ay)|^2\leq\,\| x\|_{a,\lambda}^2\| y\|_{a,\lambda}^2 ,
\end{align*}}
So by \eqref{32}, we get
$$ \mathrm{Im} (\lambda \mu \varphi(x^* ay) + (1-\lambda)\mu \overline{\varphi(ax)}\varphi(ay) )=0, $$
and hence
$$ \lambda \mu \varphi(x^* ay) + (1-\lambda)\mu \overline{\varphi(ax)}\varphi(ay ) = \| x\|_{a,\lambda}\| y\|_{a,\lambda}.$$
Therefore
$$ \lambda  \varphi(x^* ay) + (1-\lambda) \overline{\varphi(ax)}\varphi(ay ) =\overline{\mu} \| x\|_{a,\lambda}\| y\|_{a,\lambda}.  $$
\item[$ 2) \Rightarrow 1) $]
Suppose that there exist $ \varphi \in \mathcal{S}_a (\mathcal{A}) $ and $ \mu \in \mathbb{T} $ such that
\begin{align*}
 \mu(\lambda \varphi(x^* ay) + (1-\lambda) \overline{\varphi(ax)}\varphi(ay)) =\| x\|_{a,\lambda}\| y\|_{a,\lambda}.
\end{align*}
Hence
{\small\begin{align*}
\| x\|_{a,\lambda}\| y\|_{a,\lambda} &= \mu (\lambda \varphi(x^* ay) +(1-\lambda) \overline{\varphi(ax)}\varphi(ay)) \\
&=|\mu(\lambda \varphi(x^* ay)+(1-\lambda)\overline{\varphi(ax)}\varphi(ay))|\\
&\leq \, \lambda| \varphi(x^* ay)|+(1-\lambda)|\overline{\varphi(ax)}\varphi(ay)|\\
&\leq  \lambda \varphi(x^* ax)^{\frac{1}{2}} \varphi(y^* ay)^{\frac{1}{2}} + (1-\lambda ) |\varphi(ax)| \, |\varphi(ay)| \\
&\leq \sqrt{\lambda \varphi(x^* ax)+(1-\lambda)|\varphi(ax)|^2 } \sqrt{\lambda \varphi(y^* ay)+(1-\lambda)|\varphi(ay)|^2 }\\
&\leq  \sqrt{\lambda \varphi(x^* ax)+(1-\lambda)|\varphi(ax)|^2 }  \| y\|_{a,\lambda} \leq \| x\|_{a,\lambda}\| y\|_{a,\lambda} .
\end{align*}}
So we get $ \| x\|_{a,\lambda} = \sqrt{\lambda \varphi(x^* ax)+(1-\lambda)|\varphi(ax)|^2 } $.
By a similar way, we conclude that
 $ \| y\|_{a,\lambda}=\sqrt{\lambda \varphi(y^* ay)+(1-\lambda)|\varphi(ay)|^2 } $.
On the other hand,
\begin{align*}
\| x\|_{a,\lambda} \| y\|_{a,\lambda} &= \lambda\mu \varphi(x^* ay)\big)+(1-\lambda)\mu \overline{\varphi(ax)} \varphi(ay)\\
&=\mathrm{Re} (\lambda\mu \varphi(x^* ay)+(1-\lambda)\mu \overline{\varphi(ax)} \varphi(ay)).
\end{align*}
Hence
{\small\begin{align*}
&(\|x\|_{a,\lambda} +\|y\|_{a,\lambda})^2 \geq  \|x+\mu y\|_{a,\lambda}^2 \geq \lambda \varphi((x+\mu y)^* a(x+\mu y)) +(1-\lambda) |\varphi(a(x+\mu y))|^2 \\
&= \lambda (\varphi(x^* ax) + \mu \varphi(x^* ay)+\overline{\mu} \varphi(y^* ax) +\varphi(y^* ay)) \\
&+ (1-\lambda)(|\varphi(ax)|^2 + |\varphi(ay)|^2 + \mu \overline{\varphi(ax)} \varphi(ay) +\overline{\mu} \varphi(ax) \overline{\varphi(ay)}) \\
&= \big(\lambda \varphi(x^* ax) + (1-\lambda)|\varphi(ax)|^2\big) + \big(\lambda \varphi(y^* ay) +(1-\lambda) |\varphi(ay)|^2 \big) \\
&+ \lambda \big(\mu \varphi(x^* ay)+\overline{\mu} \varphi(y^* ax)\big) + (1-\lambda)  \big(\mu \overline{\varphi(ax)} \varphi(ay) +\overline{\mu} \varphi(ax) \overline{\varphi(ay)}\big) \\
 &= (\lambda \varphi(x^* ax) + (1-\lambda)|\varphi(ax)|^2 ) + (\lambda \varphi(y^* ay) +(1-\lambda) |\varphi(ay)|^2 ) \\
&+ 2 \mathrm{Re} (\lambda \, \mu \varphi(x^* ay)+(1-\lambda) \mu \, \overline{\varphi(ax)} \varphi(ay))\\
&= \| x\|_{a,\lambda}^2 + \| y\|_{a,\lambda}^2 +2 \| x\|_{a,\lambda} \| y\|_{a,\lambda} = (\| x\|_{a,\lambda}+ \| y\|_{a,\lambda})^2 .
\end{align*}}
Therefore $\|x+\mu y\|_{a,\lambda}=\|x\|_{a,\lambda}+\|y\|_{a,\lambda}$.
\end{proof}
\begin{corollary}\label{30}
Let $ \mathcal{A} $ be a unital $ C^* $-algebra, $ a \in \mathcal{A} $ be a positive and invertible element and $ x,y \in \mathcal{A} $. Then the following statements are equivalent:
\begin{itemize}
\item[1)] $ x \|_{a , \lambda } y $;
\item[2)] There exists $ \varphi \in \mathcal{S}_a (\mathcal{A}) $ such that
\begin{align*}
&(i) \lambda \varphi (x^* ay) +(1-\lambda ) \overline{\varphi(ax)} \varphi (ay)=\overline{\mu} \| x \|_{a , \lambda } \,  \| y \|_{a , \lambda} ;\\
 &(ii) \| x \|_{a , \lambda}^2=\lambda \varphi (x^* ax) +(1-\lambda)|\varphi(ax)|^2 ;\\
&(iii) \| y \|_{a , \lambda}^2=\lambda \varphi (y^* ay) +(1-\lambda )|\varphi(ay)|^2 .
\end{align*}
\end{itemize}
\end{corollary}
Taking $ \lambda=0 $ and $ \lambda=1 $ in Theorem \ref{002}, we have the following characterizations for $ a $-norm parallelism and $ a $-numerical radius parallelism in unital $ C^* $-algebras.
\begin{corollary}\label{38}
Let $ \mathcal{A} $ be a unital $ C^* $-algebra, $ a \in \mathcal{A} $ be a positive and invertible element and $ x,y \in \mathcal{A} $.
Then
\begin{itemize}
\item[1)]
$ x \|_{v_a} y $ if and only if there are $ \varphi \in \mathcal{S}_a (\mathcal{A}) $ and $\mu\in\mathbb{T}$ such that
\begin{align*}
\overline{\varphi (ax)}\varphi (a y)=\overline{\mu}v_a (x) \,  v_a (y) .
\end{align*}
\item[2)]
$ x \|_{a} y $  if and only if there are $ \varphi \in \mathcal{S}_a (\mathcal{A}) $ and $\mu\in\mathbb{T}$ such that
\begin{align*}
 \varphi (x^*ax)=\|x\|_a^2, \quad  \varphi (y^*ay)= \|y\|_a^2\quad{\rm and}\quad \varphi (x^* a y )=\overline{\mu} \| x\|_a \, \| y\|_a.
\end{align*}
\end{itemize}
\end{corollary}
%%%%%%%%%%%%%%%%%%%%%%%%%%%%%%%%%%%%%%%%%%%%%%%
%%%%%%%%%%%%%%%%%%%%%%%%%%%%%%%%%%%%%%%%%%%%%%%
%%%%%%%%%%%%%%%%%%%%%%%%%%%%%%%%%%%%%%%%%%%%%%%%%%%
\begin{theorem}
Let $ \mathcal{A} $ be a unital $ C^* $-algebra, $ a \in \mathcal{A} $ be a positive and invertible and $ x,y \in \mathcal{A}$. Then the following statements are equivalent:
\begin{itemize}
\item[1)] $ x \|_{a,\lambda} 1_{\mathcal{A}} $ $ (1_{\mathcal{A}}\|_{a,\lambda} x  )$;
\item[2)] $ x $  is $ a $-normaloid; i.e., $v_a(x)=\|x\|_a$;
\item[3)] $  x \|_{a,\lambda} x^\sharp $.
\end{itemize}
\end{theorem}
\begin{proof}
\item[$ 1) \Rightarrow 2) $]
Let $ x \|_{a,\lambda} 1_{\mathcal{A}} $. Then Theorem \ref{002} implies that there are $ \varphi \in \mathcal{S}_a (\mathcal{A}) $ and $ \mu\in \mathbb{T} $ such that
$$ \lambda \varphi(x^* a 1_{\mathcal{A}}) +(1-\lambda) \overline{\varphi(ax)} \varphi(a 1_{\mathcal{A}})=\overline{\mu} \|x\|_{a,\lambda} \|1_{\mathcal{A}}\|_{a,\lambda}= \overline{\mu} \|x\|_{a,\lambda}. $$
Hence by \eqref{008}, we get
\begin{align*}
 \|x\|_{a,\lambda}&= \mu(\lambda \varphi(x^* a) +(1-\lambda) \varphi (x^* a)) = \mu \varphi(x^* a) = |\mu \varphi(x^* a)|\\
& = |\varphi(x^* a)|= |\varphi(ax)| \leq v_a (x) \leq \|x\|_{a,\lambda} .
\end{align*}
Thus $ v_a (x)= \|x\|_{a,\lambda} $, and therefore $ v_a (x)= \|x\|_{a} $. The implication $ 1) \Rightarrow 2) $ with the assumption $ 1_{\mathcal{A}}\|_{a,\lambda} x  $ can be proved by a similar argument.
\item[$ 2) \Rightarrow 1) $]
Assume that $ x $ is $ a $-normaloid. Then there exists $ \varphi \in \mathcal{S}_a (\mathcal{A}) $ such that $ \|x\|_a = \|x\|_{a,\lambda} = v_a (x)=|\varphi(ax)| $. Moreover,
$$|\varphi(ax)|= |\overline{\varphi(ax)}|=|\varphi(x^* a)|=|\varphi(x^* a 1_\mathcal{A})| . $$
So there exists $ \mu \in \mathbb{T} $ such that $ \|x\|_{a,\lambda} =\mu \varphi(x^* a1_{\mathcal{A}}) $. Hence
\begin{align*}
\|x\|_{a,\lambda} =& \mu \varphi(x^* a 1_{\mathcal{A}})
= \mu \big( \lambda \varphi(x^* a 1_{\mathcal{A}}) +(1-\lambda) \varphi(x^* a 1_{\mathcal{A}})\varphi(a1_{\mathcal{A}}) \big)\\
& = \mu (\lambda \varphi(x^* a 1_{\mathcal{A}}) +(1-\lambda) \overline{\varphi(ax)}\varphi(a1_{\mathcal{A}})),
\end{align*}
and so Theorem \ref{002} shows that $ x\|_{a,\lambda} 1_{\mathcal{A}} $. The same method can be used to prove that $ 1_{\mathcal{A}}\|_{a,\lambda} x  $.
\item[$ 1) \Rightarrow 3) $]
Let $ x\|_{a,\lambda} 1_{\mathcal{A}} $. Then there exists $ \mu \in \mathbb{T} $ such that $ \| x+\mu 1_{\mathcal{A}}\|_{a,\lambda} = \|x\|_{a,\lambda} + 1 $. Hence there exists $ \varphi \in \mathcal{S}_a (\mathcal{A}) $ such that
{\small\begin{align*}
(\|x\|_{a,\lambda} + 1)^2 &=  \| x+\mu 1_{\mathcal{A}}\|_{a,\lambda}^2 =\lambda \varphi((x+\mu 1_{\mathcal{A}})^* a (x+\mu 1_{\mathcal{A}})) + (1-\lambda) |\varphi(a(x+\mu 1_{\mathcal{A}}))|^2 \\
&= \lambda \varphi(x^* ax) +(1-\lambda)|\varphi(ax)|^2 +2{\rm Re}(\overline{\mu} \varphi(ax))+1 \\
&\leq  \|x\|_{a,\lambda}^2 + 2 | \varphi(ax) |+1 \\
&\leq  \|x\|_{a,\lambda}^2 + 2 v_a(x) +1 \\
&\leq  \|x\|_{a,\lambda}^2 + 2\|x\|_{a,\lambda} +1 =(\|x\|_{a,\lambda} + 1)^2 .
\end{align*}}
Thus
$$ \varphi(\mu x^* a)=\|x\|_{a,\lambda}=\varphi(\overline{\mu} ax). $$
Since $ \|x\|_{a,\lambda}=\|x^\sharp\|_{a,\lambda} $, by the Cauchy-Schwartz inequality, we conclude that
\begin{align*}
\|x\|_{a,\lambda} + \|x^\sharp\|_{a,\lambda}&= \varphi(\overline{\mu} ax + \mu x^* a ) = \varphi (\overline{\mu} ax + \mu a x^\sharp) \\
 &=\varphi (a(\overline{\mu} x + \mu x^\sharp)) \\
&\leq \sqrt{\varphi (\overline{\mu} x + \mu x^\sharp)^* a (\overline{\mu} x + \mu x^\sharp)}
\leq \, \|\overline{\mu} x + \mu x^\sharp \|_{a}.
\end{align*}
But
{\small\begin{align*}
a(\overline{\mu} x + \mu x^\sharp)=\overline{\mu}a x + \mu ax^\sharp =\overline{\mu} (x^\sharp)^* a + \mu x^* a= ( \overline{\mu} (x^\sharp)^* + \mu x^* )a = (\mu x^\sharp +\overline{\mu} x)^* a.
\end{align*}}
Hence $ \mu x^\sharp +\overline{\mu} x $ is $ a $-selfadjoint, and so $ \|\mu x^\sharp +\overline{\mu}x\|_a = \|\mu x^\sharp +\overline{\mu}x\|_{a,\lambda} $. It follows that
$$ \|x\|_{a,\lambda} + \|x^\sharp\|_{a,\lambda} \leq \|\overline{\mu} x + \mu x^\sharp \|_{a,\lambda}=\| x + \mu^2 x^\sharp \|_{a,\lambda} \leq \| x\|_{a,\lambda} + \| x^\sharp \|_{a,\lambda}.  $$
Therefore $ \| x + \mu^2 x^\sharp \|_{a,\lambda} =\| x\|_{a,\lambda} + \| x^\sharp \|_{a,\lambda} $ for $ \mu^2 \in \mathbb{T} $. This yields that $ x\|_{a,\lambda} 1_{\mathcal{A}} $.
\item[$ 3) \Rightarrow 1) $]
Suppose that $ x\|_{a,\lambda} x^\sharp $, i.e., $ \| x+\mu x^\sharp \|_{a,\lambda} = \|x\|_{a,\lambda} + \|x^\sharp\|_{a,\lambda} =2\|x\|_{a,\lambda} $ for some $\mu\in\mathbb{T}$. Since
$$ a(x+\mu x^\sharp) =ax+\mu a x^\sharp = (x^\sharp)^* a + \mu x^* a =(x^\sharp +\overline{\mu}x)^* a=(\mu x^\sharp +x)^* a ,$$
we conclude that $ x+\mu x^\sharp $ is $ a $-selfadjoint. So
\begin{align*}
2\|x\|_{a,\lambda}=\| x+\mu x^\sharp\|_{a,\lambda} &=v_a (x+\mu x^\sharp) = | \varphi (a(x+\mu x^\sharp))|= |\varphi(ax)+\mu \varphi(ax^\sharp)| \\
& \leq 2 |\varphi(ax)| \leq 2 v_a (x) \leq 2\| x\|_{a,\lambda}
\end{align*}
for some $\varphi\in\mathcal{S}(\mathcal{A})$. Therefore $ \|x\|_{a,\lambda} =|\varphi(ax)| $, and so $ \|x\|_{a,\lambda}=\mu \varphi(ax) $ for some $ \mu \in \mathbb{T} $. It follows that
\begin{align*}
 \|x\|_{a,\lambda} +1 &= \varphi (\mu ax +a) = \varphi(a(\mu x+1_{\mathcal{A}}))\\
&\leq v_a (\mu x +1_{\mathcal{A}}) \leq \| \mu x +1_{\mathcal{A}}\|_{a,\lambda} = \| x +\overline{\mu}1_{\mathcal{A}}\|_{a,\lambda} \leq \|x\|_{a,\lambda} +1
\end{align*}
Thus $ \| x+\mu 1_{\mathcal{A}}\|_{a,\lambda}=\|x\|_{a,\lambda} +1 $ for some $ \mu \in \mathbb{T} $, and therefore $ x\|_{a,\lambda}1_{\mathcal{A}} $.
\end{proof}
%%%%%%%%%%%%%%%%%%%%%%%%%%%%%%%%%%%%%%%%%%%%
\begin{corollary}
Let $ \mathcal{A} $ be a unital $ C^*$-algebra, $a\in\mathcal{A}$ be a positive and invertible element. Let $ x,y \in \mathcal{A} $ are  $ a $-normaloid elements. If $ x\|_{v_a} y $, then $ x\|_{a,\lambda} y $.
\end{corollary}
\begin{proof}
Let $ x\|_{v_a} y $. Then there is some $ \mu \in \mathbb{T} $ such that $ v_a (x+\mu y) = v_a (x)+v_a(y) $. In addition, $ x,y $ are $ a $-normaloid. It follows that $ v_a (x)=\|x\|_{a,\lambda} $ and $ v_a (y)=\|y\|_{a,\lambda} $. Therefore
$$ \|x\|_{a,\lambda} + \|y\|_{a,\lambda}= v_a (x)+v_a (y)=v_a (x+\mu y) \leq \|x+\mu y\|_{a,\lambda} \leq \|x\|_{a,\lambda} + \|y\|_{a,\lambda}, $$
which imlies that $ \|x+\mu y\|_{a,\lambda}=\|x\|_{a,\lambda} + \|y\|_{a,\lambda} $ for some $ \mu \in \mathbb{T} $, and so $ x\|_{a,\lambda} y $.
\end{proof}
%%%%%%%%%%%%%%%%%%%%%%%%%%%%%%%%%%%%%%%%%%%%%%%
%%%%%%%%%%%%%%%%%%%%%%%%%%%%%%%%%%%%%%%%%%%%%%%%%%%%%%%%%%%%%%
\section{$ (a,\lambda )$-Birkhoff-James orthogonality in unital $C^*$-algebras}
\subsection{Characterization of $ (a,\lambda )$-Birkhoff-James orthogonality in unital $C^*$-algebras}

In this section we investigate the notion of 
Birkhoff-James orthogonality in unital $C^*$-algebra $\mathcal{A}$ with respect to $\|\cdot\|_{a,\lambda}$ whenever $a\in\mathcal{A}$ is positive and invertible.
\begin{definition}
Let $ \mathcal{A} $ be a unital $ C^* $-algebra,  $a\in\mathcal{A}$ be a positive element and $\lambda\in[0,1]$. An element $ x \in \mathcal{A}$ is called $ (a,\lambda )$-Birkhoff-James orthogonal to $ y \in \mathcal{A} $, denoted by $ x \perp_{a , \lambda} y $, if
$$ \| x + \xi y \|_{a, \lambda } \geq \| x\|_{a ,\lambda} \qquad \forall \xi \in \mathbb{C} . $$
\end{definition}
Note that Birkhoff-James orthogonality and numerical radius Birkhoff-James orthogonality with respect to $\|\cdot\|_a$ and $v_a(\cdot)$ is just $(a,1)$-Birkhoff-James orthogonality and $(a,0)$-Birkhoff-James orthogonality, respectively:

$\bullet$ An element $x\in\mathcal{A}$ is called  $a$-Birkhoff-James orthogonal to
another element $y\in\mathcal{A}$, denoted by $x\perp^a_{BJ}y$,  if $\|x+\xi y\|_a\geq\|x\|_a$ for all $\xi\in\mathbb{C}$.

$\bullet$ An element $x\in\mathcal{A}$ is called  $a$-numerical radius Birkhoff-James orthogonal to
another element $y\in\mathcal{A}$, denoted by $x\perp^{v_a}_{BJ}y$,  if $v_a(x+\xi y)\geq v_a(x)$ for all $\xi\in\mathbb{C}$.

Some basic proporties of $ (a,\lambda )$-Birkhoff-James orthogonality obtain in the following Proposition. The proofs are follows directly from the definition and we omit their expression.
\begin{proposition}\label{Th1}
Let $ x , y \in \mathcal{A}_a $. Then the following statements are equivalent:
\begin{itemize}
\item[1)]$x \perp_{a , \lambda} y $;
\item[2)]$ x^\sharp  \perp_{a , \lambda} y^\sharp $;
\item[3)]$ \alpha x \perp_{a , \lambda} \beta y   $ for all $ \alpha , \beta \in \mathbb{C} $; i.e., $ (a,\lambda )$-Birkhoff-James orthogonality is homogenous.

Also, note that $ (a,\lambda )$-Birkhoff-James orthogonality is non-degenerate and if $ x $ is $ (a,\lambda )$-Birkhoff-James orthogonal to $y$, then $ x $ and $ y $ are linearly independent.
\end{itemize}
\end{proposition}
%%%%%%%%%%%%%%%%%%%%%%%%%%%
%%%%%%%%%%%%%%%%%%
In the next theorem, a characterization of $(a,\lambda)$-Birkhoff–James orthogonality in unital $C^*$-algebras in terms of the elements of $\mathcal{S}_a(\mathcal{A})$ is presented.
%%%%%%%%%%
\begin{theorem}\label{33}
Let $ x,y \in \mathcal{A} $. The following statements are equivalent:
\begin{itemize}
\item[1)]$ x \perp_{a,\lambda} y $;
\item[2)] For each $  \theta \in [0,2\pi) $, there exists $ \varphi \in \mathcal{S}_a (\mathcal{A}) $ such that
\begin{align*}
\lambda \varphi (x^* a x) +(1-\lambda ) | \varphi (ax) |^2 = \| x \|_{a , \lambda}^2
\end{align*}
and
\begin{align*}
\mathrm{Re} \big(e^{i \theta } (\lambda \varphi (x^* ay) +(1-\lambda) \overline{\varphi(ax)} \varphi (ay))\big) \geq 0.
\end{align*}
\end{itemize}
\end{theorem}
\begin{proof}
\item[$ 1) \Rightarrow 2) $]
\underline{Step1.} First, we show that if
 $ \| x + \xi y \|_{a , \lambda} \geq \| x\|_{a , \lambda} $ for all $ \xi \in \mathbb{R}^+ $, then there exists $ \varphi \in \mathcal{S}_a (\mathcal{A}) $ such that
\begin{align}\label{Eq1}
\lambda \varphi (x^* a x) +(1-\lambda ) | \varphi (ax) |^2 = \| x \|_{a , \lambda}^2
\end{align}
and
\begin{align}\label{Eq2}
\mathrm{Re} ( \lambda \varphi (x^* ay) + (1-\lambda ) \overline{\varphi(ax)} \varphi (ay) ) \geq 0 .
\end{align}
If $ \| x\|_{a , \lambda} =0 $, then
$$ 0=\| x\|_{a , \lambda}^2 \, \geq \, \lambda \varphi(x^* ax) +(1-\lambda)|\varphi (ax)|^2  $$
for all $\varphi\in\mathcal{S}_a(\mathcal{A})$, and so $ \varphi(x^* ax) =|\varphi (ax)| =0 $ ( Note that for $ \lambda =0 $, we have $ |\varphi (ax)|=0 $ for all $ \varphi\in\mathcal{S}_a(\mathcal{A}) $. Then $ \varphi(x^* a)= \varphi (ax)=0 $. Hence $ \varphi(ax - x^* a)=0 $ for all $\varphi\in\mathcal{S}_a(\mathcal{A})$. It follows from \cite[Lemma 2.1]{AB} that $ ax-x^*a =0 $. Thus $ ax=x^* a $. So $ x $ is a $ a $-selfadjoint and thus $ v_a (x)=\|x\|_a $ by \eqref{008}. Then $ \|x\|_a =0 $, and so $ \varphi(x^* ax)=0 $ for all $ \varphi\in\mathcal{S}_a(\mathcal{A}) $). It follows that two side of \eqref{Eq1} and \eqref{Eq2} are equal to zero. Indeed, by the Cauchy-Schwartz inequality,
$$ | \varphi(x^* ay) | \leq \varphi(x^* ax)^{\frac{1}{2}} \varphi(y^* ay)^{\frac{1}{2}} =0,$$
and hence $ \varphi(x^* ay)=0 $.

Now, let $ \| x\|_{a , \lambda} \neq 0 $. There exists $ \varepsilon_0 \in (0,1) $ such that for each $ \varepsilon \in (0,\varepsilon_0) $, we have
\begin{align}\label{18}
\| x\|_{a , \lambda} -\varepsilon^2 \geq 0 .
\end{align}
By the assumption and \eqref{18}, for each $ \varepsilon \in (0,\varepsilon_0) $, we get
 \begin{align}\label{19}
  \| x+\varepsilon y\|_{a , \lambda} \geq \| x\|_{a , \lambda} > \| x\|_{a , \lambda} - \varepsilon^2 \geq 0 ,
 \end{align}
 and so
 \begin{align}\label{20}
  \| x+\varepsilon y\|_{a , \lambda}^2 \geq \| x\|^{2}_{a , \lambda} -2\varepsilon^2 \| x\|_{a , \lambda} + \varepsilon^4 .
 \end{align}
 Since $ \mathcal{S}_a (\mathcal{A}) $ is $ w^* $-compact, for each $ \varepsilon >0 $, there exists $ \varphi_\varepsilon \in \mathcal{S}_a (\mathcal{A}) $ such that
{\small \begin{align}\label{41}
\| x+\varepsilon y\|_{a,\lambda}^2 &= \lambda \varphi_\varepsilon ( (x+\varepsilon y)^* a(x+\varepsilon y)) + (1-\lambda) | \varphi_\varepsilon (a(x+\varepsilon y))|^2 \\
 &= \lambda \big(\varphi_\varepsilon (x^* ax) + 2\varepsilon \mathrm{Re} \varphi_\varepsilon (x^* ay) + \varepsilon^2 \varphi_\varepsilon (y^* ay) \big) \nonumber\\
 &+ (1-\lambda) \big(|\varphi_\varepsilon (ax)|^2 + \varepsilon^2 |\varphi_\varepsilon (y^* ay)|^2 + 2 \varepsilon \mathrm{Re} (\overline{\varphi_\varepsilon (ax)}\varphi_\varepsilon(ay))\big) \nonumber\\
&\leq  \big(\lambda \varphi_\varepsilon (x^* ax) + (1-\lambda) |\varphi_\varepsilon (ax)|^2 \big) + \big(\lambda \varepsilon^2 \varphi_\varepsilon (y^* ay) )+(1-\lambda) \varepsilon^2 |\varphi_\varepsilon (y^* ay)|^2\big)\nonumber\\
&+  2 \varepsilon \, \lambda \, \mathrm{Re} \, \varphi_\varepsilon (x^* ay)  +
(1-\lambda)\, \mathrm{Re}\, (\overline{\varphi_\varepsilon (ax)}\varphi_\varepsilon(ay))\nonumber \\
&\leq  \| x\|_{a,\lambda}^2 + \|y\|_{a,\lambda}^2 +  2 \varepsilon \big( \lambda  \mathrm{Re} \varphi_\varepsilon (x^* ay)  +
(1-\lambda) \mathrm{Re} (\overline{\varphi_\varepsilon (ax)}\varphi_\varepsilon(ay)\big) \label{40} \\
&\leq  \| x\|_{a,\lambda}^2 + \|y\|_{a,\lambda}^2 +
  2 \varepsilon \big(\lambda |\varphi_\varepsilon (x^* ay)| + (1-\lambda) |\varphi_\varepsilon (ax)|\, |\varphi_\varepsilon(ay)|\big) \nonumber\\
&\leq  \| x\|_{a,\lambda}^2 + \|y\|_{a,\lambda}^2 +
 2 \varepsilon \big(\lambda \varphi_\varepsilon(x^* ax)^{\frac{1}{2}}\varphi_\varepsilon(y^* ay)^{\frac{1}{2}} + (1-\lambda)  |\varphi_\varepsilon (ax)|\, |\varphi_\varepsilon(ay)|\big) \nonumber\\
  &\quad \text{(by the Cauchy-Schwartz inequality)} \nonumber \\
& \leq  \| x\|_{a,\lambda}^2 + \|y\|_{a,\lambda}^2 \\
 &+ 2\varepsilon \sqrt{\lambda \varphi_\varepsilon (x^* ax) +(1-\lambda )|\varphi_\varepsilon (ax)|^2 }\, \sqrt{ \lambda\varphi_\varepsilon (y^* ay) +(1-\lambda) |\varphi_\varepsilon (ay)|^2}\nonumber \\
&\quad \text{(by the Cauchy-Bunyakovsky-Schwartz inequality)}\nonumber\\
& \leq  \|x\|_{a,\lambda}^2 + \varepsilon^2  \|y\|_{a,\lambda}^2 + 2\varepsilon  \|x\|_{a,\lambda}^2 \|y\|_{a,\lambda}^2 =( \|x\|_{a,\lambda}^2 + \varepsilon  \|y\|_{a,\lambda}^2)^2 .\nonumber
 \end{align}}
 So
 \begin{align}\label{22}
\| x+\varepsilon y\|_{a, \lambda}^2  \leq ( \|x\|_{a,\lambda}^2 + \varepsilon  \|y\|_{a,\lambda}^2)^2.
 \end{align}
Also, $ w^* $-compactness of $ \mathcal{S}_a (\mathcal{A})$ implies that there exists $ \varphi \in \mathcal{S}_a (\mathcal{A}) $ such that $ \varphi_\varepsilon \overset{w^*}{\longrightarrow}\varphi $, whenever $ \varepsilon\rightarrow 0^+ $.
From \eqref{19} and \eqref{22} we get
\begin{align*}
\| x\|_{a, \lambda}^2 &\leq \| x+\varepsilon y\|_{a, \lambda}^2 = \lambda \varphi_\varepsilon ( (x+\varepsilon y)^* a(x+\varepsilon y)) + (1-\lambda) | \varphi_\varepsilon (a(x+\varepsilon y))|^2 \\
&\leq ( \|x\|_{a,\lambda}^2 + \varepsilon  \|y\|_{a,\lambda}^2).
\end{align*}
Then
$$ \| x\|_{a, \lambda}^2 = \lambda \varphi (x^* ax) + (1-\lambda) | \varphi (ax)|^2 . $$
Moreover, by \eqref{41} and \eqref{40} we get
{\small\begin{align*}
\| x+\varepsilon y\|_{a, \lambda}^2 =& \lambda \varphi_\varepsilon ( (x+\varepsilon y)^* a(x+\varepsilon y)) + (1-\lambda) | \varphi_\varepsilon (a(x+\varepsilon y))|^2 \\
\leq & \|x\|_{a,\lambda}^2  + \varepsilon^2 \|y\|_{a,\lambda}^2  + 2 \varepsilon\, \big(\lambda \, \mathrm{Re} \, \varphi_\varepsilon(x^* ay) +(1-\lambda) \, \mathrm{Re} (\overline{\varphi_\varepsilon(ax)}\varphi_\varepsilon(ay)\big) .
\end{align*}}
It follows from \eqref{20} that
 {\small\begin{align*}
 2 \varepsilon \big(\lambda \, \mathrm{Re} \, \varphi_\varepsilon(x^* ay) +(1-\lambda) \, \mathrm{Re} (\overline{\varphi_\varepsilon(ax)}\varphi_\varepsilon(ay) \big) \geq -2\varepsilon^2 \|x\|_{a,\lambda} + \varepsilon^4 +\varepsilon^2 \|y\|_{a,\lambda}^2 ,
\end{align*}}
 and so
{\small\begin{align*}
 \lambda \, \mathrm{Re} \, \varphi_\varepsilon(x^* ay) +(1-\lambda) \, \mathrm{Re} (\overline{\varphi_\varepsilon(ax)}\varphi_\varepsilon(ay) \geq
  -2\varepsilon \|x\|_{a,\lambda} + \dfrac{\varepsilon^3}{2} +\dfrac{\varepsilon}{2} \|y\|_{a,\lambda}^2 .
\end{align*}}
Now, by letting $ \varepsilon\rightarrow 0^+ $, we conclude that
$$ \lambda \, \mathrm{Re} \, \varphi (x^* ay) +(1-\lambda) \, \mathrm{Re} (\overline{\varphi(ax)}\varphi (ay) \geq 0 . $$

\underline{Step2.}
Assume that for each $ \xi \in \mathbb{C} $, we have
$ \| x+\xi y\|_{a,\lambda}^2 \geq \| x\|_{a,\lambda}^2 $. Fix $ \theta \in [0,2\pi) $. Take $ y_\theta =e^{i\theta}y $. Since $ \| x+|\xi| y_\theta \|_{a,\lambda}^2 \geq \| x\|_{a,\lambda}^2 $, we conclude that
$$ \mathrm{Re}\big( e^{i\theta} (\lambda\,  \varphi (x^* ay) +(1-\lambda) \overline{\varphi(ax)} \varphi(ay)\big) \geq 0 , $$
by step 1.
\item[$ 2) \Rightarrow 1) $]
\underline{Step1.} Assume that there exists $ \varphi \in \mathcal{S}_a (\mathcal{A}) $ such that
\begin{align*}
\lambda \varphi (x^* a x) +(1-\lambda ) | \varphi (ax) |^2 = \| x \|_{a , \lambda}^2
\end{align*}
and
\begin{align*}
\mathrm{Re} ( \lambda \varphi (x^* ay) + (1-\lambda ) \overline{\varphi(ax)} \varphi (ay) ) \geq 0 .
\end{align*}
Then for each $ \xi \in \mathbb{R}^+ $,
{\small\begin{align*}
\| x+\xi y\|_{a,\lambda}^2& \geq \, \lambda \varphi \big( (x+\xi y)^* a (x+\xi y)\big) +(1-\lambda)|\varphi(a(x+\xi y))|^2 \\
&=\, (\lambda \varphi(x^* ax) +(1-\lambda )|\varphi(ax)|^2 ) + \xi^2 (\lambda \varphi(y^* ay) +(1-\lambda )|\varphi(ay)|^2 ) \\
&+ 2\xi \, \mathrm{Re} (\lambda \varphi(x^* ay) +(1-\lambda ) \overline{\varphi(ax)} \varphi(ay)) \\
&\geq \, \lambda \varphi(x^* ax) +(1-\lambda )|\varphi(ax)|^2 = \| x\|_{a,\lambda}^2 .
\end{align*}}
\underline{Step2.} Assume that for each $ \theta \in [0,2\pi) $, there exists $ \varphi \in \mathcal{S}_a (\mathcal{A}) $ such that
\begin{align*}
\lambda \varphi (x^* a x) +(1-\lambda ) | \varphi (ax) |^2 = \| x \|_{a , \lambda}^2
\end{align*}
and
\begin{align*}
{\rm Re} (e^{i \theta } \lambda \varphi (x^* ay) +   e^{i \theta }(1-\lambda ) \overline{\varphi(ax)} \varphi (ay) ) \geq 0.
\end{align*}
 By the polar decomposition, for each $ \xi \in \mathbb{C} $, there exists $ \theta \in [0,2\pi ) $ such that $ \xi=|\xi| e^{i \theta} $.
Hence
{\small\begin{align*}
\| x+\xi y\|_{a,\lambda}^2& \geq \, \lambda \varphi((x+\xi y)^* a (x+\xi y)) +(1-\lambda) |\varphi(a(x+\xi y))|^2 \\
&=\, \lambda \varphi(x^* ax) + (1-\lambda)|\varphi(ax)|^2 + |\xi e^{i \theta}|^2 ( \lambda \varphi(y^* ay) +(1-\lambda)|\varphi(ay)|^2) \\
&+ 2|\xi| \mathrm{Re}(e^{i \theta} \lambda \varphi (x^* ay)+e^{i \theta} (1-\lambda) \overline{\varphi(ax)} \varphi(ay)) \\
&\geq  \lambda \varphi(x^* ax)  + (1-\lambda)|\varphi(ax)|^2 = \| x\|_{a,\lambda}^2 .
\end{align*}}
\end{proof}
%%%%%%%%%%%%%%%%%%%%%%%%%%%%%%%%%%%%%%%%%%%%%%%%%%%%%%%%
The following results are direct consequences of the previous theorem which give some characterizations of  $a$-Birkhoff–James orthogonality and $a$-numerical radius Birkhoff–James orthogonality in unital $C^*$-algebras.
\begin{corollary}\label{36}
Let $ \mathcal{A} $ be a unital $ C^* $-algebra, $ a \in \mathcal{A} $ be a positive and invertible and $ x,y \in \mathcal{A} $.
Then
\begin{itemize}
\item[1)]
$ x \perp_{BJ}^{v_a} y $ if and only if for each $\theta\in [0,2\pi)$ then there exists $ \varphi  \in \mathcal{S}_a (\mathcal{A}) $  such that
\begin{align*}
| \varphi (ax )| =v_a (x), \quad Re ( e^{i \theta } \overline{\varphi (ax )} \varphi (ay )) \geq 0 .
\end{align*}
\item[2)]
$ x \perp_{BJ}^{a} y $ if and only if for each $\theta\in [0,2\pi)$ there exists $ \varphi_  \in \mathcal{S}_a (\mathcal{A}) $  such that
\begin{align*} \varphi (x^* ax) = \| x\|_a^{2} , \quad Re ( e^{i \theta} \varphi (x^* ay)) \geq 0 .
\end{align*}
\end{itemize}
\end{corollary}
\begin{corollary}\label{C3-5}
Let $ x , y \in \mathcal{A} $ be $a$-positive elements. Then the following statements are equivalent:
\begin{itemize}
\item[1)] $ x \perp_{BJ}^{v_a} y $;
\item[(2)] There exists $ \varphi\in\mathcal{S}_a (\mathcal{A}) $ such that $ v_a (x) =\varphi (ax) $ and $ \varphi (ay)=0 $.
\end{itemize}
\end{corollary}
\begin{proof}
\item[$ 1) \Rightarrow 2) $] Let $ x \perp_{BJ}^{v_a} y $. Then  $ x \perp_{BJ}^{v_a} (-y) $, by Theorem  \ref{Th1}. So, there exists $ \varphi \in \mathcal{S}_a (\mathcal{A}) $ such that $ v_a (x) = | \varphi (ax)| $ and $ Re ( \overline{\varphi (ax)} \varphi(a(-y)) \geq 0 $. Since $ x$ is $a$-positive, we have $ | \varphi (ax )| = \varphi (ax )=v_a(x) $ and $ \varphi ( ay ) \geq 0 $. Hence
$$ 0 \leq \varphi ( ay) = \dfrac{- Re  ( \overline{\varphi (ax)} \varphi(a(-y)) }{\varphi ( ax )}  \leq 0 ,$$
and so $ \varphi ( ay) = 0 $.
\item[$ 2) \Rightarrow 1) $] Let $ \varphi\in\mathcal{S}_a (\mathcal{A}) $  be such that $ v_a (x) =\varphi (ax) $ and $ \varphi (ay)=0 $. Then for each $ \lambda \in \mathbb{C} $, we get
$$ v_a ( x + \lambda y) \geq \, | \varphi (a(x+ \lambda y ))| = | \varphi (ax) + \lambda \varphi (ay)| = | \varphi (ax )| = \varphi (ax ) = v_a (x) . $$
\end{proof}
In the following theorem we obtain the connection between $(a,\lambda)$-norm parallelism
and $(a,\lambda)$-Birkhoff-James orthogonality of the elements in $C^*$-algebras.
\begin{theorem}\label{35}
Let $ \mathcal{A} $ be a unital $ C^* $-algebra, $ a \in \mathcal{A} $ be a positive and invertible and $ x,y \in \mathcal{A}$. If $ x \|_{a,\lambda} y $, Then the following statements hold:
\begin{itemize}
\item[1)] There exists $ \mu \in \mathbb{T} $ such that
$ x \perp_{a,\lambda} (\|y\|_{a,\lambda} x - \mu \|x\|_{a,\lambda} y) $.
\item[2)] There exists $ \mu \in \mathbb{T} $ such that
$ y \perp_{a,\lambda} (\|y\|_{a,\lambda} x - \overline{\mu} \|x\|_{a,\lambda} y) $.
\end{itemize}
\end{theorem}
\begin{proof}
Let us just prove (1). The proof of (2) is similar.
\begin{itemize}
\item[1)]
Assume that $ x \|_{a,\lambda} y $. By Corollary \ref{30}, there are $ \varphi \in \mathcal{S}_a (\mathcal{A}) $ and $\mu\in\mathbb{T}$ such that
$$ \|x\|_{a,\lambda}^2 =\lambda \varphi(x^* ax)+(1-\lambda)|\varphi(ax)|^2 $$
and
$$ \lambda \varphi (x^* ay) +(1-\lambda ) \overline{\varphi(ax)} \varphi (ay)=\overline{\mu} \| x \|_{a , \lambda } \,  \| y \|_{a , \lambda}. $$
Hence
{\small\begin{align*}
&\lambda \varphi (x^* a (\|y\|_{a,\lambda} x - \mu \|x\|_{a,\lambda} y)) + (1-\lambda) \overline{\varphi(ax)} \varphi(a(\|y\|_{a,\lambda} x - \mu \|x\|_{a,\lambda} y)) \\
=&\lambda\|y\|_{a,\lambda} \varphi (x^* a x) - \lambda\mu \|x\|_{a,\lambda} \varphi (x^* ay)) +(1-\lambda) (\|y\|_{a,\lambda} |\varphi(ax)|^2 - \mu \|x\|_{a,\lambda} \overline{\varphi(ax)} \varphi(ay) ) \\
=& \|y\|_{a,\lambda} (\lambda \varphi (x^* a x) +(1-\lambda) |\varphi(ax)|^2 ) - \mu \|x\|_{a,\lambda} (\lambda \varphi (x^* ay)+(1-\lambda)\overline{\varphi(ax)} \varphi(ay) ) \\
=& \|y\|_{a,\lambda} \|x\|_{a,\lambda}^2 - \|x\|_{a,\lambda}^2  \|y\|_{a,\lambda}=0.
\end{align*}}
Thus for each $\theta\in[0,2\i)$, we have
{\small\begin{align*}
&\mathrm{Re} \Big( e^{i \theta}\big(\lambda \varphi (x^* a (\|y\|_{a,\lambda} x - \mu \|x\|_{a,\lambda} y)) + (1-\lambda) \overline{\varphi(ax)} \varphi(a(\|y\|_{a,\lambda} x - \mu \|x\|_{a,\lambda} y)) \big) \Big)=0.
\end{align*}}
Consequently, Theorem \ref{33} complete the proof.
\end{itemize}
\end{proof}
%%%%%%%%%%%%%%%%%%%%%%%%%%%%%%%%%%%%%%%%%%%%%%%%%%%
\noindent\textbf{Acknowledgements}

The authors would like to express their sincere gratitude to the anonymous referee for his/her helpful comments.

%%%%%%%%%%%%%%%%%%%%%%%%%%%%%%%%%%%%%%%%%%%
%%%%%%%%%%%%%%%%%%%%%%%%%%%%%%%%%%%%%%%%%%%%%%%%%%%%%%%%%
\bibliographystyle{amsplain}

\end{document}